\documentclass{amsart}
\usepackage[utf8]{inputenc} 
\usepackage[T1]{fontenc}    
\usepackage[backref=page]{hyperref}       
\usepackage{amsfonts}       
\usepackage{nicefrac}       
\usepackage{microtype}      
\usepackage{enumerate}

\usepackage{tikz}

\hypersetup{
  colorlinks   = true, 
  urlcolor     = teal, 
  linkcolor    = teal, 
  citecolor   = teal 
}

\usepackage[color=teal!70!white, bordercolor=teal, textsize=footnotesize,obeyFinal]{todonotes}

\usepackage{marvosym}
\usepackage{amsmath}
\usepackage{amsgen}
\usepackage{amsbsy}
\usepackage{amsopn}
\usepackage{amsthm}
\usepackage{amscd}
\usepackage{amsfonts}
\usepackage{amssymb}
\usepackage{mathrsfs}
\usepackage[capitalize]{cleveref}

\usepackage{ifsym}

\usepackage{tikz-cd}
\usepackage{epigraph}

\usepackage[usestackEOL]{stackengine}
\renewcommand{\epigraphsize}{\small}
\newcommand{\mytextformat}{\epigraphsize\itshape}
\newcommand{\mysourceformat}{\epigraphsize\scshape}

\let\originalepigraph\epigraph 
\renewcommand\epigraph[2]{%
  \setbox0=\hbox{\stackon{\textit{\mytextformat\Longstack{#1}}}%
    {\mysourceformat\scshape\Longstack{#2}}}%
  \ifdim\wd0>.8\linewidth\wd0=.8\linewidth\fi%
  \setlength{\epigraphwidth}{\wd0}%
  \originalepigraph{\textit{#1}}{\textsc{#2}}%
}

\usepackage{color}

\usepackage{graphicx}
\usepackage{comment}

\usepackage{latexsym}

\newtheorem{Theorem}{Theorem}
\newtheorem{Lemma}[Theorem]{Lemma}

\newtheorem{Corollary}[Theorem]{Corollary}

\newtheorem{Proposition}[Theorem]{Proposition}

\newtheorem{Definition}{Definition}

\newtheorem{Remark}{Remark}

\newcommand{\R}{\bf R}
\def \C{{\mathbf C}}

\def\Sl{\widetilde{\rm SL}_{2}}

\newcommand{\cout}[1]{}
\def\co{\colon\thinspace}

\begin{document}

\title[Asymptotic dimension and geometric decompositions]{Asymptotic dimension and geometric \\ decompositions in dimensions 3 and 4}

\author[H. Contreras Peruyero]{H. Contreras Peruyero }
\author[P. Su\'arez-Serrato]{P. Su\'arez-Serrato \\  \today}
\address{Instituto de Matem\'aticas, Universidad Nacional Aut\'onoma de M\'exico UNAM, Mexico Tenochtitlan}
\email{pablo@im.unam.mx} 

\address{Centro de Ciencias Matem\'aticas, Universidad Nacional Aut\'onoma de M\'exico UNAM, Antigua Carretera a P\'atzcuaro \# 8701, Col. Ex Hacienda San Jos\'e de la Huerta, 58089, Morelia, Michoac\'an, Mexico}

\email{haydeeperuyero@matmor.unam.mx}

\thanks{HCP thanks the support of a DGAPA-UNAM postdoctoral research grant 2023}

\begin{abstract}
We show that the fundamental groups of smooth $4$-manifolds that admit geometric decompositions in the sense of Thurston have asymptotic dimension at most four, and equal to 4 when aspherical.
We also show that closed $3$-manifold groups have asymptotic dimension at most 3. 
 Our proof method yields that the asymptotic dimension of closed $3$-dimensional Alexandrov spaces is at most 3.
We thus obtain that the Novikov conjecture holds for closed $4$-manifolds with such a geometric decomposition and closed  $3$-dimensional Alexandrov spaces. 
Consequences of these results include a vanishing result for the  Yamabe invariant of certain $0$-surgered geometric $4$-manifolds and the existence of zero in the spectrum of aspherical smooth $4$-manifolds with a geometric decomposition.
\end{abstract} 

\maketitle

\section{Introduction}

The uniformization theorem for topological surfaces, proved by Koebe \cite{Koebe} and Poincar\'e \cite{Poincare}, showed that geometric structures can effectively classify distinct families of manifolds. 
Their geometric classification scheme involves the three constant sectional curvature $2$-dimensional geometries.

In dimension three, constant sectional curvature manifolds are insufficient to include all $3$-manifolds. 
W.T. Thurston defined a model geometry as a complete, simply connected Riemannian manifold $\mathbb{X}$ such that the group of isometries acts transitively on $\mathbb{X}$ and contains a discrete subgroup with a finite volume quotient.
A manifold $X$ is said to be {\it geometrizable}, in the sense of Thurston, if $X$ is diffeomorphic to a connected sum of manifolds which admit a decomposition into pieces, each modeled on a Thurston geometry.

In dimension three there are eight \emph{model geometries}, and manifolds modeled in these serve as building blocks that when assembled produce a global description of a $3$-manifold.
Comprehensive descriptions of the model geometries in dimension $3$ may be found in Thurston's book \cite{Thurston}, and in a survey by Scott \cite{Scott}.
The success of the geometrization program in dimension three by Thurston and Hamilton--Perelman \cite{Perelman} leads us to wonder about the nature of geometrizable manifolds in higher dimensions.
Filipkiewicz \cite{Fi} classified all maximal four dimensional model geometries.
The following list includes all of the four-dimensional Thurston geometries that admit finite volume quotients: 
\[   
\begin{array}{ccccccc}	
     \mathbb{S}^{4}, & \mathbb{C}{\rm P}^{2}, & \mathbb{S}^{3}\times \mathbb{E}, & \mathbb{H}^{3}\times \mathbb{E},   &
 \Sl \times \mathbb{E}, & \mathbb{N}il^{3}\times \mathbb{E}, \\ \mathbb{N}il^{4}, &

 \mathbb{S}^{2} \times\mathbb{E}^{2}, & \mathbb{H}^{2}\times\mathbb{E}^{2},
& \mathbb{S}ol^{4}_{m,n}, & \mathbb{S}ol^{4}_{1},& \mathbb{S}ol^{4}_{0},  \\
 \mathbb{S}^{2}\times \mathbb{S}^{2},  & \mathbb{S}^{2}\times\mathbb{H}^{2}, & \mathbb{E}^{4}, &  \mathbb{F}^4, &
 \mathbb{H}^{4}, & \mathbb{H}^{2}\times\mathbb{H}^{2}, & \mathbb{H}_{\C}^{2}  
\end{array}
\]

\begin{figure}[h!]
\centering
\usetikzlibrary{mindmap}
\scalebox{0.70}{
\begin{tikzpicture}
  \path[mindmap,concept color=black,text=white]
    node[concept] (c1) {$ \mathbb{H}^{4}$}
    [clockwise from=0]
    
    child[concept color=green!55!gray] {
      node[concept] {$\mathbb{H}^{3}\times \mathbb{E}$}
      [clockwise from=-30]
      child { node[concept] { $\mathbb{H}^{2}\times\mathbb{E}^{2}$} }
      child { node[concept] {$\Sl \times \mathbb{E}$} }
    }
    
    child[concept color=blue!50!gray] {
      node[concept] {$\mathbb{H}_{\C}^{2}$}
      [clockwise from=-30]
      child { node[concept] {$\mathbb{F}^4$} }
      child { node[concept] {$\mathbb{F}^4$} }
    }
    child[concept color=teal!65!white] { node[concept] {$\mathbb{S}^{3}\times \mathbb{E}$} }
    child[concept color=violet!65!gray] { node[concept] {$\mathbb{C}{\rm P}^{2}$} };
\end{tikzpicture}
}
\caption{A sketch of a $4$-manifold $X$ whose connected summands admit a decomposition into Thurston geometries.
Each circular region represents a geometric manifold, tagged with its model geometry. 
Colors indicate parts that are either geometric, or have a proper geometric decomposition. 
The purple region represents a submanifold $X_1$, that decomposes into pieces modeled on the geometries $\mathbb{H}_{\C}^{2}$ and $\mathbb{F}^4$, glued along nilpotent boundaries.
In the green region we see a submanifold $X_2$ that decomposes into $\mathbb{H}^{3}\times \mathbb{E}$, $\mathbb{H}^{2}\times\mathbb{E}^{2}$, and $\Sl \times \mathbb{E}$ pieces, 
 glued along flat boundaries.
The strips joining different colored areas represent connected sums. Let the central, real hyperbolic piece (in black) be $X_0$. Then $X= X_0 \# X_1 \#X_2 \# S^3\times S^1 \# \mathbb{C}P^2$.
} \label{fig:Geomz-4mfd}
\end{figure}
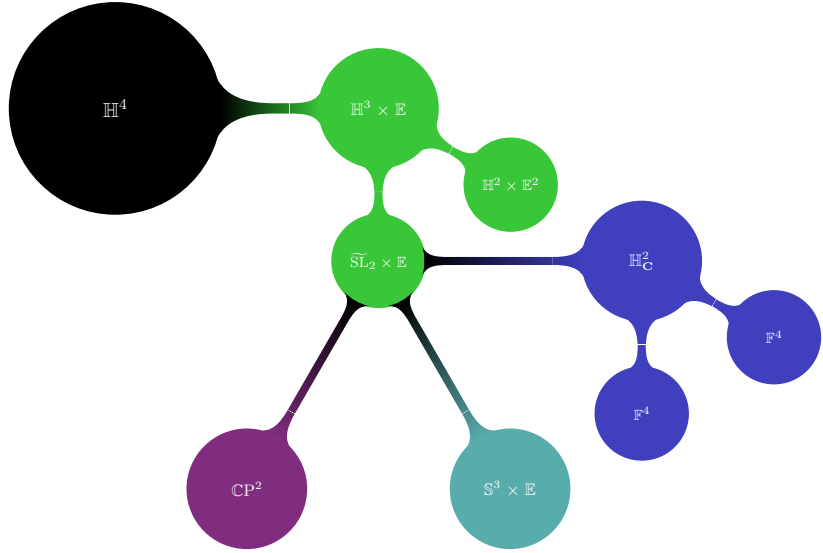

Detailed explanations and examples for all these geometries are available in the work of Hillman \cite[p.133]{Hil} and Wall \cite{Wall}. 
To keep our exposition short, we recommend interested readers consult those sources.
In \cref{fig:Geomz-4mfd} below we show a schematic example of one of these manifolds. 

We have previously studied the minimal volume entropy problem \cite{SS09}, and the existence of Einstein metrics \cite{CPSS22} on manifolds in this family. 

 Gromov defined the asymptotic dimension, ${\rm asdim}\ \Gamma$, of a metric space $\Gamma$ as a coarse analogue of the Lebesgue covering dimension \cite{Gro2} (for details see \S \ref{sec:def-asdim} below).  
 
 Here, we show:

\begin{Theorem}\label{thm:4Thurstongeom-finite-asdim}
Let $X$ be one of the following;
\begin{enumerate}[i)] 
\item a closed orientable $4$-manifold that is geometric in the sense of Thurston;
\item a closed orientable $4$-manifold that is geometrizable in the sense of Thurston. 
\end{enumerate}
Then, the asymptotic dimension of $\pi_1(X)$ is at most 4, and it equals 4 when $X$ is aspherical.
\end{Theorem}

In the case of $3$-manifolds, by the same methods, we can show the following:

\begin{Theorem}\label{thm:3mfds-asdim-leq-3}
    Let $Y$ be a closed $3$-manifold, then  ${\rm asdim}\ \pi_1(Y) \leq 3$, and it equals $3$ when $Y$ is aspherical.
\end{Theorem}

Our proof uses the geometrization of $3$-manifolds, specifically the description of $\pi_1(Y)$ as a graph of groups.
Even though it may be known to some experts, \cref{thm:3mfds-asdim-leq-3} improves upon the available published bounds for ${\rm asdim}\ \pi_1(Y)$ \cite{MacSis, EnMa}. 
In related work, Ren showed that such a $\pi_1(Y)$ has finite decomposition complexity \cite{Ren}, however a bound was not made explicit.

Alexandrov spaces are a generalization of smooth manifolds with bounded curvature, in the sense that they include all limits of sequences of smooth manifolds with sectional curvatures bounded below. 
Briefly, they are locally complete, locally compact, connected length spaces, that satisfy a lower curvature bound in the triangle-comparison sense.

 Following the same terminology as for $3$-manifolds, an Alexandrov $3$-dimensional space $Y$ is called geometric, with a given model Thurston geometry, if $Y$ can be written as a quotient of that geometry by some cocompact lattice. 
 A closed Alexandrov $3$-dimensional space is said to admit a geometric decomposition, if there exists a collection of spheres, projective planes, tori and Klein bottles that decompose $Y$ into geometric pieces.
 A geometrization theorem for Alexandrov $3$-spaces was shown by F. Galaz-Garc\'ia and Guijarro.
 They showed that a closed three-dimensional Alexandrov space admits a geometric decomposition into geometric three-dimensional Alexandrov spaces \cite[Theorem 1.6]{GGG}.  
 Moreover, they proved that an Alexandrov $3$-space $Y$ may be presented as the quotient of a smooth $3$-manifold $Y^*$ by the action of an isometric involution \cite[Lemma 1.8]{GGG}.
 
 The universal cover of an Alexandrov space is by definition the simply connected cover with the induced metric structure making the covering map into a local isometry.
 The fundamental group of a compact Alexandrov space here will be seen as a discrete group of isometries of its universal cover. 
 More details and related rigidity results for Alexandrov $3$-spaces may be found in the recent survey by Núñez-Zimbrón \cite{Núñez-Zimbrón}.
 
We obtain the following result as a consequence of \cref{thm:3mfds-asdim-leq-3}:
\begin{Theorem}\label{thm:3Alex-asdim-bound-3}
Let $Y$ be a closed $3$-dimensional Alexandrov space, then  $${\rm asdim}\ \pi_1(Y) \leq 3.$$   
\end{Theorem}

Our proof crucially uses the quasi-isometric invariance property of ${\rm asdim}$, applied to a specific action of a group on the universal covering space of the $3$-manifold $Y^*$ that produces $Y$ after quotienting out by the isometric involution. 
    
Recall that a space is aspherical if its universal cover is contractible. 

 Our main theorems yield information about the Baum--Connes, Novikov, and zero in the spectrum conjectures, contained in the next Corollaries. 
 
While we give more details about these topics in \S 2 below, we point interested readers to available reviews of these conjectures, for example, by Yu \cite{Yu}, Davis \cite{DavisAspects}, Ferry--Ranicki--Rosenberg \cite{FRR}, and Weinberger \cite{Weinberger}. 

\begin{Corollary}\label{cor:consequences-fin-asdim}
Let $X$ be a manifold from \cref{thm:4Thurstongeom-finite-asdim}, then:
\begin{enumerate}[i)]
    \item The coarse Baum--Connes conjecture holds for $X$.
    \item The Novikov conjecture holds for $X$.
    \item If $X$ is aspherical then its universal cover $\widetilde{X}$ has a zero in the spectrum.
\end{enumerate}   
\end{Corollary}

The second item above includes manifolds constructed using complex hyperbolic pieces, which are not included in previous related work on higher graph manifolds \cite{FLS, CS19, BJS17}. 
Moreover, these earlier results do not apply to dimension four, because they depend on arguments from surgery theory that are not known to hold yet for groups of exponential growth \cite{FT}.

Lott showed that for a closed geometric 4-manifold $X$, zero is in the Laplace--Beltrami spectrum of $\widetilde{X}$ \cite[Proposition 18]{Lott}. 
By comparison, \cref{cor:consequences-fin-asdim} is shown by different methods, it subsumes previous work on aspherical geometric manifolds, and further includes all the aspherical manifolds in \cref{thm:4Thurstongeom-finite-asdim}.
 
 Similarly, \cref{thm:3Alex-asdim-bound-3} has the following consequences:

\begin{Corollary}\label{cor:consequences-fin-asdim-3Alex}
Let $Y$ be a be a closed $3$-dimensional Alexandrov space, then:
\begin{enumerate}[i)]
    \item The coarse Baum--Connes conjecture holds for $Y$.
    \item The Novikov conjecture holds for $Y$.
    \item If $X$ is aspherical then its universal cover $\widetilde{X}$ has a zero in the spectrum.
\end{enumerate}   
\end{Corollary}

Previous related work showing the Novikov conjecture holds for singular spaces includes that of Ji on buildings \cite{Ji06}, and on  torsion free arithmetic subgroups of connected, rational, linear algebraic groups \cite{Ji}.

A conjecture attributed to Gromov--Lawson--Rosengberg states that there do not exist Riemannian metrics with positive scalar curvature on compact aspherical manifolds. 

 As a consequence of \cref{thm:4Thurstongeom-finite-asdim} and work of Yu \cite{Y} (and, alternatively, of Dranishnikov \cite{D03}) the aspherical manifolds in Theorem \ref{thm:4Thurstongeom-finite-asdim} do not admit Riemannian metrics of positive scalar curvature. 
 This result was recently shown for all  aspherical smooth 4-manifolds by Chodosh--Li \cite{CL23}, and by Gromov \cite{Gro3}. 
 Nevertheless, our methods provide an independent proof for the manifolds in \cref{thm:4Thurstongeom-finite-asdim}.

A natural approach to understanding how topology and geometry are coupled is by minimizing the curvature that a Riemannian manifold may have.
One way to achieve this is to minimize a norm of a curvature tensor.
Let $(M,g)$ be a compact Riemannian manifold with a smooth metric $g$. 
Consider a conformal class of Riemannian metrics, 
 $   \gamma:= [g] = \{u\cdot g \ | \ M\overset{u}{\longrightarrow} \R^+\}$.

The \emph{Yamabe constant} of $(M,g)$ is defined as, 
 \begin{equation*} 
 \mathcal{Y}(M, \gamma):= \underset{g\in \gamma}{\inf} \frac{\int_M {\rm Scal}_g d{\rm vol}_g}{({\rm Vol} (M, g))^{2/n}}.
 \end{equation*}
Here ${\rm Scal}_g$ denotes the scalar curvature and $d{\rm vol}_g$ the volume form of $g$. 
The Yamabe invariant of a manifold $M$
is then defined to be
 
 $\mathcal{Y}(M):= \underset{\gamma}{\sup} \ \mathcal{Y}(M, \gamma)$. 

We next present an application of \cref{thm:4Thurstongeom-finite-asdim} to the study of the Yamabe invariant.

\begin{Corollary}\label{cor:zero-Yamabe}
  Let $X$ be a manifold from \cref{thm:4Thurstongeom-finite-asdim}. 
  If the geometric pieces of $X$ are modelled on the geometries,
    \[
    \begin{array}{ccccc}
     \mathbb{E}^{4}, & \mathbb{H}^{3}\times \mathbb{E},  & \mathbb{H}^{2}\times\mathbb{E}^{2},  & \mathbb{N}il^{4}, &  \mathbb{S}ol^{4}_{1},  \\
    \Sl \times \mathbb{E}, & \mathbb{N}il^{3}\times \mathbb{E}, & \mathbb{S}ol^{4}_{m,n}, & \mathbb{S}ol^{4}_{0}, & or\, \,  \mathbb{F}^4,   
    \end{array}    
    \]
     then the Yamabe invariant of $X \# k(S^{3}\times S^1)$, with $k\in \{ 0, 1, 2, \ldots \, \}$, vanishes. 
\end{Corollary}

This improves upon a result by the second named author, \cite[Proposition 2.6 (i)]{SS08}, covering only closed $\mathbb{E}^{4}, \mathbb{H}^{3}\times \mathbb{E},$ or $ \mathbb{H}^{2}\times\mathbb{E}^{2}$ manifolds. 
 In those cases the existence of a nonpositive sectional curvature metric obstructs the existence of positive scalar curvature metrics.
 
 Wall \cite{Wall} showed there is a close relationship between geometric structures and complex surfaces.
 So in \cref{cor:zero-Yamabe} there is some overlap with the work of LeBrun, who, as a part of a {\it tour-de-force} of results on the Yamabe invariant, showed that compact complex surfaces of Kodaira dimension $0$ or $1$ have null Yamabe invariant \cite{LeBrun}.
 For compact complex surfaces that admit a geometric structure listed in \cref{cor:zero-Yamabe}, we now have an independent proof that their Yamabe invariant is zero.
 For example, the compact complex surfaces known as Inoue surfaces are exactly those admitting one of the geometries $\mathbb{S}ol^{4}_{0}$ or $\mathbb{S}ol^{4}_{1}$ \cite{Wall}. Albanese recently showed that Inoue surfaces  have zero Yamabe invariant \cite{Alba}, and \cref{cor:zero-Yamabe} now gives an alternative proof. 

Finally, we take this opportunity to include the following result that recovers part of LeBrun's theorem mentioned above \cite{LeBrun}, and for which we can now produce a simple proof (given what is needed for the previous Corollary).

\begin{Lemma}\label{lem:ComplexSurfaces-ZeroYam}
Let $X$ be an aspherical compact complex surface of Kodaira dimension at most 1 and which is not of class VII.
Then $\mathcal{Y}(X)=0$.
\end{Lemma}

Here, as usual, the Kodaira dimension $\kappa = \limsup_{m\to \infty} (\log(P_{m}(X)) / \log m)$, with $P_{m}(X)$ the dimension of the space of holomorphic sections of the $m$-th tensor power of the
canonical line bundle of $X$, and $\kappa := -\infty$ if $P_{m}(X)=0$ for all $m$.
Surfaces of Kodaira dimension $-\infty$ that are not K\"ahler are called surfaces of class VII. 
These include Inoue surfaces  with vanishing second Betti number (featured in \cref{cor:zero-Yamabe}), Hopf surfaces (known to be geometric, but are not aspherical), certain compact elliptic surfaces, and surfaces with a {\it global spherical shell}, which have positive second Betti number and are not aspherical.
These are conjecturally all the minimal surfaces of class VII.

The relevant definitions for the concepts appearing in \cref{cor:consequences-fin-asdim}, \cref{cor:consequences-fin-asdim-3Alex}, and \cref{cor:zero-Yamabe} are found in \cref{subsec:BaumConnes}, \cref{subsec:Novikov}, \cref{subsec:Zero}, and  \cref{subsec:Yamabe}. 
The proofs of each of the items in \cref{cor:consequences-fin-asdim} and \cref{cor:consequences-fin-asdim-3Alex} appear as \cref{cor:Coarse-BC-conj}, \cref{cor:Novikov-Conj}, and \cref{cor:zerospec}. 

The proofs of \cref{thm:4Thurstongeom-finite-asdim} and \cref{thm:3mfds-asdim-leq-3} rely on close examinations of the fundamental groups involved, both are in \cref{sec:main-proofs}. 
We use various properties and operations on groups to bound the asymptotic dimension from both sides.
Related results are available for higher graph manifolds, due to the second named author in collaborations with Connell \cite{CS19}, and with B\'arcenas and Juan Pineda \cite{BJS17}. 
All these strategies are reminiscent of the original work of Wall \cite[\S 12]{Wall2} on codimension 1 splittings along a hypersurface, and of Cappell \cite{Cappell, Cappell2} on amalgamated products.
However, those arguments from classical surgery theory need to assume the dimension of the manifold is at least 5.

\section*{Acknowledgements} We thank Daniel Juan Pineda for clarifying the hypothesis needed to apply Yu's result. We also thank Chris Connell for pointing out how to rephrase the statement of one of our Corollaries.

\section{Preliminaries and proofs of Corollaries \ref{cor:consequences-fin-asdim}, \ref{cor:consequences-fin-asdim-3Alex}, and \ref{cor:zero-Yamabe}}  

\subsection{Definition of asymptotic dimension}\label{sec:def-asdim}
Gromov introduced the concept of asymptotic dimension of a metric space $(X , d)$ \cite{Gro2}. There are several equivalent definitions, we state the following: 

\begin{Definition}\label{def:asdim}
 We say that the asymptotic dimension of $(X,d)$ does not exceed 
 $n$, written ${\rm asdim}\ X \leq n$, if for each $D > 0$, there exist $B \geq 0$ and families $\mathcal{U}_0 ,\ldots , \mathcal{U}_n$ of subsets which form a cover of $X$ such that:
 \begin{enumerate}[i)]
     \item for all $i\leq n$ and all $U$ in $\mathcal{U}_i$, their diameter satisfies ${\rm diam}\ (U)\leq B$;
     \item for all $i\leq n$ and all $U$ and $V$ in $\mathcal{U}_i$, if $U\neq V$ then $d(U, V) > D$.
 \end{enumerate}
 \end{Definition}

\begin{figure}[h!]
\centering
\begin{tikzpicture}
\usetikzlibrary{patterns}
\draw [pattern color=teal,pattern=bricks] (0,0) circle (2cm); 
\end{tikzpicture}
\caption{This covering by bricks is helpful in understanding why the asdim of ${\bf R}^2$ is at most $2$. Observe that a point in the plane either lies in the interior of a brick, or it lies on a the boundary of a brick. In the former case there is a neighborhood contained in the brick. In the latter, there are two option, it either lies exactly at the point where 3 bricks meet or not. In either of these cases, the neighborhood of the point will intersect at most $(\dim ({\bf R}^2) + 1 )$ bricks.} \label{fig:asdim-bricks}
\end{figure}
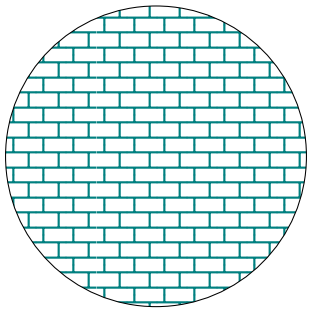

While internalizing this definition may take some time, we recommend consulting the friendly and accessible exposition by Bell \cite{Bell-asdimOH}.
In \cref{fig:asdim-bricks} we see a specific cover, by bricks on a plane, illustrating the two points of \cref{def:asdim}. 
First, all bricks are isometric, so their diameter is the same. 
Second, different bricks need to be translated at least a distance $D$ to match, and this quantity depends on the size of the brick (itself determined by their diameter $B$).

Let $\Gamma$ be a finitely generated group and let $S$ be a finite generating set. 
The word length with respect to $S$, denoted by $l_{S}$, of an element $\gamma\in\Gamma$ is the smallest integer $n\geq 0$ for which there exist $s_{1},...,s_{n}\in S\cup S^{-1}$ such that $\gamma=s_1...s_n$. 
The word metric, denoted by ${\rm d}_S$ is defined as ${\rm d}_{S}(\gamma_{1},\gamma_{2})=l_{S}(\gamma_{1}^{-1}\gamma_{2})$.
A finitely presented group, equipped with the word metric, is a metric space.
We refer the reader to the work of Bell and Dranishnikov \cite{BD} for multiple examples of groups and spaces with finite asymptotic dimension.
For a finitely generated group $\Gamma$, the asymptotic dimension is a group property, i.e., it is independent of the choice of generators \cite[Corollary 51]{BD}.

\begin{Lemma}\label{lem:finite-asdim}
    \begin{enumerate}[i)]
        \item \cite[Example 9.6]{Roe03} The euclidean $n$-dimensional space $\mathbb{E}^n$  has asymptotic dimension equal to $n$.
        \item \cite{Roe05} The real hyperbolic $n$-dimensional space $\mathbb{H}^n$ has asymptotic dimension equal to $n$. 
        \item \cite[Proposition 60]{BD} Let $\Gamma$ be a finitely generated group. Then ${\rm asdim}\ \Gamma = 0$ if and only if $\Gamma$ is finite.
    \end{enumerate}
\end{Lemma}

A pair of metric spaces $(X_1 , d_1 ), (X_2,d_2)$ are called {\it quasi-isometric}, if there exists a map $f : X_1 \to X_2$ and a constants $B> 0$ and $C \geq 1$ such that,
\begin{enumerate}
    \item for every pair of points $x, y$ in $X_1$,
    \[
\frac{1}{B}\cdot d_1 (x,y) - C \leq d_2 (f(x),f(y)) \leq B\cdot d_1 (x,y) + C,   
    \] 
    \item  every point of $X_2$ lies  within a $C$-neighborhood of the image $f(X_1)$.
\end{enumerate}

A well known property of the asymptotic dimension is that it is an invariant of the quasi-isometry type of a finitely generated group $\Gamma$ \cite{Gro2}. 
As a consequence of Milnor--\v{S}varc Lemma, if $M$ is a compact Riemannian manifold with universal cover $\widetilde{M}$ and finitely generated group $\pi_1(M)$, then 
$\widetilde{M}$ is quasi-isometric to $\pi_1(M)$ (with the word metric).
Therefore \cite[Corollary 56]{BD},
\begin{equation}\label{eqn:universal-asdim}
    {\rm asdim}\left( \widetilde{M} \right) ={\rm asdim}(\pi_1(M)). 
\end{equation}

As an illustrative example, and because we will need it later on, let us focus now on the case of surface groups and show---the well known fact---that they have asymptotic dimension bounded above by $2$:

\begin{Lemma}\label{lem:asdim-2D}
        Let $\Gamma$ be the fundamental group of a closed $2$-manifold, then  ${\rm asdim}\ \Gamma \leq 2$.
    \end{Lemma}

    \begin{proof}        
        By the uniformization theorem for surfaces, $\Gamma$ may be represented as either the trivial group, a flat $2$-manifold group $\Gamma_{F}$, or the fundamental group of a genus $g\geq 2$ surface with a hyperbolic metric $\Gamma_{H}$. 
        Observe that $\Gamma_{F}$ is quasi-isometric to $\mathbb{E}^2$, and $\Gamma_{H}$ is quasi-isometric to $\mathbb{H}^2$.
        Therefore, by \cref{lem:finite-asdim} and \cref{eqn:universal-asdim} above, in all these cases  ${\rm asdim}\ \Gamma \leq 2$.       
    \end{proof}

The following result was shown by Carlsson and Goldfarb: 
\begin{Lemma}\cite[Corollary 3.6]{CG04-2}\label{lem:asdim-LieLattice}
    Let $\Gamma$ be a compact lattice in a connected Lie group $G$, and $K$ be its maximal compact subgroup. Then ${\rm asdim} (\Gamma) =  {\rm dim} (G/K).$
\end{Lemma}

An action of a discrete group $\Gamma$ on a metric space $X$ is proper if for every compact subset $B\subset X$ and for all but finitely many $\gamma$ in $\Gamma$, the intersection $\gamma(B)\cap B = \emptyset$. 
Let $\Gamma$ be a discrete group acting properly on a proper metric space $X$. 
Then, the asymptotic dimensions of $\Gamma$ and $X$ satisfy the following relationship.

\begin{Theorem}\cite[Proposition 2.3]{Ji}\label{thm:proper-action}
Let $(M,d)$ be a proper metric space. If a finitely generated group $\Gamma$ acts properly and isometrically on $M$, then for any point $x\in M$, the map $(\Gamma, \text{d}_S)\longrightarrow (\Gamma x, \text{d})$, $\gamma \longrightarrow \gamma\cdot x$ is a coarse equivalence, and hence 
\begin{equation*}
    {\rm asdim}\ \Gamma \leq {\rm asdim}\ M.
\end{equation*} 
\end{Theorem}

The following extension theorem of Bell and Dranishnikov \cite{BD} covers the case of an exact sequence:

\begin{Theorem}\cite[Theorem 63]{BD}\label{thm:extension}
    Let $1 \longrightarrow K \longrightarrow G \longrightarrow H\longrightarrow 1$ be an exact sequence with $G$ finitely generated. Then
    \begin{equation}\label{eq:extension}
        {\rm asdim}\ G \ \leq {\rm asdim}\ H \ + \ {\rm asdim}\ K.
    \end{equation} 
\end{Theorem}

The previous theorem is crucial for geometric decompositions that are injective at the level of the fundamental group, as these determine a splitting into a graph of groups.
There are also some cases, involving decompositions into $\mathbb{H}^{2}\times \mathbb{H}^{2}$ pieces, that fail to be $\pi_{1}$-injective. 
We will treat that situation with the following result.

\begin{Theorem}\cite[Finite Union Theorem]{BD2}\label{thm:finite-union}
Suppose that a metric space is presented as a union of subspaces $A\cup B$, then 
\begin{equation}\label{eq:finite-union}
    {\rm asdim}\ A \cup B \leq \max \{ {\rm asdim}\ A ,\ {\rm asdim}\ B \}.
\end{equation}   
\end{Theorem}

\subsection{Asymptotic dimension of products}

Next we recall the definition of a coarse space, following Roe's book \cite{Roe03}.

\begin{Definition}
    Let $X$ be a set. A collection of subsets $\mathcal{E}$ of $X\times X$  is called a coarse structure, and the elements of $\mathcal{E}$ are called entourages, if the following axioms are satisfied:
    \begin{enumerate}[i)]
        \item A subset of an entourage is an entourage.
        \item A finite union of entourages is an entourage.
        \item The diagonal $\Delta_{X}:= \{ (x,x) \ | \ x\in X \} $ is an entourage.
        \item The inverse $E^{-1}$ of an entourage $E$ is an entourage:
        \[E^{-1} : = \left\{ (y,x)\in X\times X  \ | \ (x,y) \in E \right\} \]
        \item The composition $E_{1}E_{2}$ of entourages $E_1$ and $E_2$ is an entourage:
        \[ E_{1}E_{2} := \left\{ (x,z)\in X\times X \, | \, \exists \ {y\in X},\, (x,y)\in E_{1},\, \text{and} \,\, (y,z)\in E_{2} \right\}  \]
    \end{enumerate}
    The pair $(X,\mathcal{E})$ is called a coarse space.
\end{Definition}

For example, topological manifolds $M$ are coarse spaces, where entourages may be defined as neighborhoods of points in $M\times M$ \cite[Chapter 2]{Roe03}.

The following was shown by Grave \cite[Proposition 20]{Gra}:

\begin{Proposition}\label{prop:asdim-products}
    Let $(X, \mathcal{E}_{X})$ and $(Y, \mathcal{E}_{Y} )$ be
coarse spaces. Then
\begin{equation}\label{eqn:asdim-products}
    {\rm asdim}(X \times Y, \mathcal{E}_{X} * \mathcal{E}_{Y} ) \leq {\rm asdim}(X, \mathcal{E}_{X}) + {\rm asdim}(Y, \mathcal{E}_{Y} ).
\end{equation}
\end{Proposition}

However, in general the equality in equation (\ref{eqn:asdim-products}) does not hold (see \cite{Gra, BD}). 

\subsection{Fundamental groups of geometrizable manifolds}\label{subsec:graph-gps}

Let $M$ be an orientable smooth four manifold which admits a proper geometric decomposition. 
A standard argument using the Seifert--van Kampen theorem shows $\pi_1(M)$ is isomorphic to an amalgamated product $A\ast_{C}B$ or to an  HNN-extension $A\ast_{C}$. 

Here $A$ is the fundamental group of one of the geometric pieces. 

Let $\Gamma$ be a graph with vertex set $V$ and directed edge set $E$.
A graph of groups over $\Gamma$ is an object $\mathcal{G}$ that assigns to each vertex $v$ a group $G_v$, and to each edge $e$ a group $G_e$, together with two injective homomorphism $\phi_e:G_e\longrightarrow G_{i(e)}$ and $\phi_{\bar{e}}:G_e\longrightarrow G_{t(e)}$. 
Here $\bar{e}$ is the edge with reverse orientation, the vertex $i(e)$ is the initial vertex of $e$ and the vertex $t(e)$ is the final vertex of $e$. 

An orientable smooth four-manifold which admits a proper, $\pi_1$-injective, geometric decomposition has a fundamental group that is isomorphic to a graph of groups constructed as an iterated amalgamated product \cite{Hil, SS09}. 

Observe that the only smooth $4$-manifolds that admit geometric decompositions which are not $\pi_1$-injective are those with irreducible $\mathbb{H}^{2}\times\mathbb{H}^{2}$ pieces \cite{Hil}. 

Bell and Dranishnikov proved the following results about the asymptotic dimensions of amalgams \cite[Theorem 82]{BD}.

\begin{Theorem}\label{thm:amalgams-fin-asdim}
  Let $A$ and $B$ be finitely generated groups and let $C$ be a subgroup of both, then: 
     \begin{equation}\label{eqn:amalgams-fin-asdim}
     {\rm asdim}(A\ast_{C}B)\leq \max\{ {\rm asdim}\ A, {\rm asdim}\ B, {\rm asdim}\ C+1 
     \}    
     \end{equation}
\end{Theorem}

\subsection{Aspherical geometrizable $4$-manifolds}\label{sec:aspherical}

Hillman obtained the following classification of closed aspherical 4-manifolds with a geometric decomposition:

\begin{Theorem}\cite[Theorem 7.2]{Hil} \label{thm:aspherical-manifolds}
    If a closed 4-manifold $M$ admits a geometric decomposition then either:
    \begin{enumerate}
        \item $M$ is geometric; or
        \item $M$ is the total space of an orbifold with general fibre $S^2$ over a hyperbolic 2-orbifold; or
        \item the components of $M\backslash \cup S$ all have geometry $\mathbb{H}^{2}\times\mathbb{H}^{2}$; or
        \item the components of $M\backslash \cup S$ have geometry $\mathbb{H}^{4},\ \mathbb{H}^{3}\times \mathbb{E}^{1}, \ \mathbb{H}^{2}\times \mathbb{E}^{2}$ or $\Sl \times \mathbb{E}^{2}$; or
        \item the components of $M\backslash \cup S$ have geometry $\mathbb{H}_{\C}^{2}$ or $\mathbb{F}^{4}$.
    \end{enumerate}
    In cases (3), (4) or (5) $\chi(M)\geq 0$ and in cases (4) or (5) $M$ is aspherical.
\end{Theorem}
In the geometric case (1) $M$ is aspherical only when its model geometry is aspherical, thus it must be modeled on $ \mathbb{E}^{4},\ \mathbb{H}^{4}  ,\ \mathbb{H}^{3}\times \mathbb{E},\ \mathbb{H}^{2}\times \mathbb{E}^{2}, \ \mathbb{H}^{2}\times \mathbb{H}^{2}, \ \mathbb{H}^{2}_{\C}, \Sl \times \mathbb{E}, \ \mathbb{N}il^{3}\times \mathbb{E},  \ \mathbb{N}il^{4}, \  \mathbb{S}ol^{4}_{1}, \ \mathbb{S}ol^{4}_{m,n},$ or $ \ \mathbb{S}ol^{4}_{0}$. 
In case (2) $M$ is never an aspherical manifold. 
Hillman expresses precise conditions under which such an orbifold bundle with a geometric decomposition is not geometric \cite[Theorem 10.2]{Hil}.

In case (3) we may or may not obtain aspherical manifolds, although some easy constructions are well known to produce aspherical examples \cite{Hil}.

\subsection{Hyperbolicity and relative hyperbolicity}\label{sec:hyp-rel-hyp}

A geodesic metric space $(X, {\rm d} )$ is called $\delta$-hyperbolic for $\delta\geq 0$ if ${\rm d}(x',y')\leq \delta$ whenever $x,y,z\in X$, $x'$ and $y'$ lie on the geodesics from $z$ to $x$ and $y$ respectively, and the following inequality is satisfied
\[ {\rm d}(x',z) = \ {\rm d}(y',z) \leq (1/2)\left( {\rm d}(x,z) + \ {\rm d}(y,z) - \ {\rm d}(x,y)\right). \]

Let $\Gamma$ be a finite group with finite generating set $S$.
The Cayley graph of $\Gamma$, with respect to $S$, is the graph $C(\Gamma,S)$ whose vertices are the elements of $\Gamma$, and whose edge set is $\left\{ (\gamma,\gamma\cdot s)|\gamma\in \Gamma, s\in S\backslash\{e\}\right\}$.
We say that $\Gamma$ is hyperbolic if the Cayley graph $C(\Gamma,S)$ associated to $\Gamma$ is a $\delta$-hyperbolic metric space, for some $\delta >0$. 
Gromov observed that hyperbolic groups have finite asymptotic dimension \cite{Gro2}, a short proof was made available by Roe \cite{Roe05}.

\begin{Theorem}\label{thm:hyp-fin-asdim}
    Finitely generated hyperbolic groups have finite asymptotic dimension.
\end{Theorem}

Let $\Gamma$ be a group and consider a collection of subgroups $\{H_{\lambda}\}_{\lambda\in \Lambda}$ of $\Gamma$, indexed by $\Lambda$. 
Let $X$ be a subset of $\Gamma$, we say that $X$ is a relative generating set of $\Gamma$ with respect to the collection $\{H_{\lambda}\}_{\lambda\in \Lambda}$ if $\Gamma$ is generated by $X\cup \left(\cup_{\lambda} H_{\lambda}\right)$. 
Let $F(X)$ be the free group with basis $X$. Then the group $\Gamma$ can be expressed as the quotient group of the free group $F = \left( \ast_{\lambda\in\Lambda} H_\lambda\right) \ast F(X)$. 
We say that the group $\Gamma$ has a relative presentation 
$\langle\, X, H_{\lambda}, \lambda\in \Lambda \,\,|\,\, R=1, R\in\mathcal{R}\,\rangle$
if the kernel $N$ of the natural homomorphism $\epsilon: F\longrightarrow \Gamma$ is a normal closure of a subset $\mathcal{R} \in N$ in the group $F$. 

If $X$ and $\mathcal{R}$ are finite then we say that the group $\Gamma$ is finitely presented relative to the collection of subgroups $\{H_{\lambda}\}_{\lambda\in \Lambda}$. 

Let $\mathcal{H}=\bigsqcup_{\lambda\in\Lambda} ( H_{\lambda} \backslash \{1\})$. 
A word $W$ in the alphabet $X\cup \mathcal{H}$ that represents $1$ in the group $\Gamma$ admits an expression in terms of the elements of $\mathcal{R}$ and $F$ as follows,
\begin{equation}\label{eq:representationW}
    W=_{F} \prod_{i=1}^{k}f_{i}^{-1}R_{i}^{\pm 1}f_{i}.
\end{equation}

Here $=_{F}$ denotes the equality in the group $F$, $R_{i}\in\mathcal{R}$, and $f_{i}\in F$ for $i=1,...,k$. 
The relative area of $W$, denoted by ${\rm Area}^{rel}(W)$ is the smallest number $k$ in a representation of the form in \cref{eq:representationW}.

A group $\Gamma$ is hyperbolic relative to a collection of subgroups $\{H_{\lambda}\}_{\lambda\in \Lambda}$ if $\Gamma$ is finitely presented relative to the collection and there is a constant $L>0$ such that for any word $W\in X\bigcup \left(\bigsqcup_{\lambda} ( H_{\lambda} \backslash \{1\})\right) $ that represent the identity in $\Gamma$ then ${\rm Area}^{rel}(W)\leq L ||W||$.  

The next result we need was shown by Dahmani--Yaman \cite[Corollary 0.2]{DY} for groups that are hyperbolic relative to a family of virtually nilpotent subgroups, and by Osin \cite[Theorem 1.2]{O} in a more general form.

\begin{Theorem}\label{thm:rel-hyp-fin-asdim}
    Let $\Gamma$ be a finitely generated group that is hyperbolic relative to a finite collection of subgroups $\{H_{\lambda}\}_{\lambda\in \Lambda}$. If each of the groups $H_{\lambda}$ has finite asymptotic dimension then ${\rm asdim}\ \Gamma<\infty$.
\end{Theorem} 

\subsection{Nagata dimension}\label{sec:Nagata-dim}

Let $X$ be a metric space and consider a family $\mathcal{B}=(B_i)_{i\in I}$ of subsets of $X$, with index set $I$. 
For some constant $D\geq 0$, the family $\mathcal{B}$ will be called $D$-bounded if, for all $i\in  I$, ${\rm diam}\ B_{i}:= \sup\{ {\rm d}(x,x') : x,x'\in B_{i}\} \leq D$. 

The multiplicity of the family is defined as the infimum over all integers $n\geq 0$ such that every point in the metric space $X$ is in at most $n$ elements of $\mathcal{B}$. 
Let $s>0$ be a constant, the $s$-multiplicity of the the family $\mathcal{B}$ is the infimum over all $n$ such that every subset of $X$ of diameter $\leq s$ intersects at most $n$ elements of $\mathcal{B}$.

\begin{Definition}
    Let $X$ be a metric space. 
    The Nagata dimension of $X$, denoted by $\dim_N X$, is the infimum of all integers $n$ such that  there exists a constant $c$ such that for all $s>0$, $X$ has a ($c\cdot s$)-bounded  covering with $s$-multiplicity at most $n+1$. 
\end{Definition}

In Figure \ref{fig:asdim-bricks}, we have a family $\mathcal{B}=(B_i)_{i\in I}$ of bricks of sides $l_1 \leq l_2$. 
These bricks constitute a $(c\cdot s)$-bounded collection of subsets of $\mathbf{R}^2$, where $c$ equals $\sqrt{(l{_1}^2+l_{2}^{2})}/l_{1}$.
Now, to see that the multiplicity of that family is at most 3 we need to observe the following cases. 
Let $x$ be a point in $\mathbf{R}^{2}$ inside a brick and consider a ball $B(x,r)$ centered on $x$ of radius $r \leq  l_1/2$. Then $B(x,r)$ intersects just one of the bricks. 
Now, let $x$ be a point in the boundary of two or tree bricks and consider again a ball $B'(x,r)$ of radius $r \leq l_1/2$ with center in $x$. 
Then $B'(x,r)$ intersects two bricks if $x$ lies on the boundary of exactly two bricks, and it intersects three bricks if $x$ lies on the corner of a brick. 
We have exhibited a family of subsets of $X$ with multiplicity at most 3, and therefore the Nagata dimension of $\mathbf{R}^2$ is at most 2.

\begin{Remark}
    By the definition, the Nagata dimension is an upper bound for the asymptotic dimension, 
    \begin{equation}\label{eqn:Nagata-asdim-bound}
        \dim_N X \geq {\rm asdim}\, X.
    \end{equation}
\end{Remark}

The following result by Lang--Schlichenmaier \cite{Lang} will be useful later on.

\begin{Theorem}\cite[Theorem 3.7]{Lang}\label{thm:Hadamard-Nagata-dim}
    Let $X$ be an $n$-dimensional Hadamard manifold whose sectional curvature $K$ satisfies $-b^2 \leq K \leq -a ^{2}$ for some positive constants $b\geq a$. 
    Then $\dim_N X = n$.
\end{Theorem}

For example, we obtain the following:

\begin{Corollary}\label{cor:Nagata-complexhyperbolic}
     The Nagata dimension of $\mathbb{H}_{\C}^{2}$ equals 4.
\end{Corollary}

\begin{proof}
    Recall that the sectional curvature of the Bergman metric on $\mathbb{H}_{\C}^{2}$ is bounded between $-4$ and $-1$.  
    As the real dimension of $\mathbb{H}_{\C}^{2}$ equals $4$, then \cref{thm:Hadamard-Nagata-dim} yields $\dim_N \mathbb{H}_{\C}^{2} =4$.
\end{proof}

\subsection{Lower cohomological bounds}

There are several equivalent definitions of cohomological dimension. Consider the following, based on K.S. Brown's book \cite[\S\ VIII]{Brown}, where the definition of group cohomology $H^{n}(\Gamma,\mathbf{Z})$ may also be found.

\begin{Definition}
    The cohomological dimension of a group $\Gamma$, denoted by ${\rm cd}(\Gamma)$, is defined as
    \begin{equation*}
        {\rm cd}(\Gamma) = \sup\left\{ n: H^{n}(\Gamma,\mathbf{Z})\neq 0 \right\}. 
    \end{equation*}
\end{Definition}

The dimension of an aspherical manifold provides an upper bound for the cohomological dimension of its fundamental group.

\begin{Proposition}\cite[Proposition 8.1]{Brown}
    Suppose that $Y$ is a $d$-dimensional $K(\Gamma,1)$-manifold (possibly with boundary). Then
    \begin{enumerate}
        \item ${\rm cd}(\Gamma) \leq d$, with equality if and only if $Y$ is closed (i.e., compact and without boundary).
        \item If $Y$ is compact then $\Gamma$ has a finite classifying space $B\Gamma$. 
    \end{enumerate}
\end{Proposition}

As a consequence of this proposition, if $M$ is an aspherical manifold, then 
\begin{equation}\label{eq:aspherical-cd-equals-dim}
{\rm cd}(\pi_{1}(M))=\dim M.    
\end{equation}

Dranishnikov showed the following:

\begin{Proposition}\cite[Proposition 5.10]{D09}\label{prop:FP-asdim-cd}
    Let $\Gamma$ be a finitely presented discrete group such that its classifying space $B\Gamma$ is dominated by a finite complex. Then 
   \begin{equation}\label{eq:asdim-cd-bound} 
    {\rm asdim\ \Gamma} \geq {\rm cd}(\Gamma).
    \end{equation}
\end{Proposition}

In the case of aspherical manifolds we obtain:

\begin{Lemma}\label{lmm:asph-cohdim}
   Let $M$ be an aspherical manifold. Then 
   \begin{equation}\label{eqn:asdim-cd}
       {\rm asdim}\ \pi_{1}(M) \geq {\rm cd} (\pi_{1}(M)) = \dim M.
   \end{equation}
\end{Lemma}

\begin{proof}
This is mentioned by Gromov in his asymptotic invariants of infinite groups essay \cite[pg. 33]{Gro2}.
A proof also follows by observing that, for a fundamental group $\pi_{1}(M)$ of a compact aspherical manifold, $M$ itself is a finite model for the classifying space $B\pi_1(M)$.
Therefore, by Proposition \ref{prop:FP-asdim-cd}, inequality \eqref{eqn:asdim-cd} holds.
\end{proof}

\subsection{Properties of Alexandrov spaces of dimension $3$}\label{subsec:Alex-3spaces}

A metric space $(X,d)$ is called a {\it length space} if for every $x,y\in X, \, [d(x,y)=\inf\{L(\gamma): \gamma(a)=x,\ \gamma(b)=y \}$. 
Here the infimum is taken over all continuous curves $\gamma:[a,b]\rightarrow X$, and $L(\gamma)$ denotes the length of the curve $\gamma$, defined as,
\begin{equation*}
L(\gamma) = \sup_{F} \left\{ \sum_{i=1}^{n-1}d(\gamma(t_{i}),\ \gamma(t_{i+1}))\right \} 
\end{equation*} 
\noindent where the supremum runs over all finite partitions $F$ of $[a,b]$.

Observe that if the length metric space $(X,d)$ is complete and locally compact then there exists at least one geodesic between each pair of points $x,y\in X$. 

Let $k$ be a real number. 
We will call a complete, simply-connected $2$-dimensional Riemannian manifold of constant curvature $k$ a {\it model space}, and denote it by $M_{k}^{2}$. 
Depending on the sign of $k$ the space $M_{k}^{2}$ is isometric to one of the following \cite{Koebe, Poincare}:

\begin{enumerate}
\item If $k>0$, we have a sphere of constant curvature $k$, $\mathbb{S}^{2}_{k}$.
\item If $k=0$, we have the Euclidean plane of null curvature, $\mathbb{E}_{k}^{2}$.
\item If $k<0$, we have a hyperbolic plane of constant curvature $k$, $\mathbb{H}_{k}^{2}$.
\end{enumerate}

Let $|\, ,|$ be the usual length metric on the corresponding model space.
Consider a geodesic triangle $pqr$ in the length space $(X,d)$. 
That is, $pqr$ is a collection of three points $p,q,r \in X$ and the segments connecting them, $[pq],[qr]$ and $[rp]$, are geodesics. 
Given a geodesic triangle $pqr$ in $X$, the geodesic triangle $\bar{p},\bar{q},\bar{r}$ in the model space $M_{k}^{2}$ is a comparison triangle for $pqr$ if $d(p,q)=|p,q|$, $d(q,r)=|q,r|$, and $d(r,p)=|r,p|$. 

A length space $(X,d)$ is said to have curvature bounded by below by $k\in\mathbb{R}$ if, for every $x\in X$, there exists an open neighborhood $U\subset X$ of $x$ such that for every geodesic triangle $pqr$ in $X$ and every comparison triangle $\bar{p},\bar{q},\bar{r}$ in the model space $M_k^{2}$, for all $s\in[p,q]$ and $\bar{s}\in[\bar(p),\bar(q)]$ such that $d(p,s)=|\bar{p},\bar{s}|$, then we have that $d(r,s)\geq |\bar{r}-\bar{s}|$. 

An Alexandrov space is a complete and locally compact length space $(X,d)$ with curvature bounded below by some $k\in\mathbb{R}$.

  The following geometrization theorem for Alexandrov $3$-spaces was shown by F. Galaz-Garc\'ia and Guijarro.

  \begin{Theorem}[\cite{GGG}]\label{thm:geomz-Alex}
    A closed three-dimensional Alexandrov space admits a geometric decomposition into geometric three-dimensional Alexandrov spaces.  
  \end{Theorem} 

  Moreover, they showed that an Alexandrov $3$-space $Y$ may be presented as the quotient of a smooth $3$-manifold $Y^*$, under the action of an isometric involution $\iota:Y^*\to Y^*$ \cite[Lemma 1.8]{GGG}. The fixed points of $\iota$ descend under the quotient to the singular points $\mathcal{S}(Y)$ of the Alexandrov structure on $Y$. 

We will need the following Lemma, included here for completeness.

\begin{Lemma}\label{lem:3Alex-Univ-Cover}
The universal coverings  $\widetilde{Y^*}$ and $\widetilde{Y}$, of the spaces $Y^*$ and $Y$ mentioned immediately above (and in the same order), are the same. This universal covering space is unique up to covering isomorphism.
\end{Lemma}

Although this is most likely evident to experts, we include a brief proof.
\begin{proof}
    In fact, the space $Y$, being presented as a global quotient of the $3$-manifold $Y^*$ is a very good orbifold. 
    In this light, the quotient map $Y^*\to Y$ is an orbifold covering map.
    Now observe the fact, first recorded by Thurston in his notes \cite{Thurston-notes} (with further details also provided by Choi \cite[Proposition 8]{Choi}), that there exists a universal covering orbifold, and that it is unique up to covering isomorphism. 
    A standard argument proves uniqueness, and extensive details are also available  \cite[Proposition 9]{Choi}.
\end{proof}

\subsection{Coarse Baum--Connes conjecture}\label{subsec:BaumConnes}
Let $M$ be a manifold and $\Gamma=\pi_1(M)$.
Recall that a metric space is called proper if closed, bounded sets are compact.
The group $\Gamma$ endowed with the word metric is a proper metric space. Consider the $C^{\ast}$--algebra $C^{\ast}\Gamma$. The coarse assembly map is defined as 
\begin{equation*}
    \mu_{X} : \ KX_{\ast}(\Gamma) \longrightarrow K_{\ast}(C^{\ast}(\Gamma))
\end{equation*}
where $K_{\ast}(C^{\ast}(\Gamma))$ denotes the $K$-theory of the $C^{\ast}$-algebra and $KX_{\ast}(\Gamma)$ is the limit of the $K$-homology groups (see \cite{V}).
A metric space is said to have bounded geometry if for every $r>0$ the cardinality of balls of radius $r$ is uniformly bounded.
The coarse Baum--Connes conjecture states that if a proper metric space has bounded geometry then the coarse assembly map is an isomorphism. 
Yu \cite{Y} proved the following:

\begin{Theorem}\label{thm:CBC}\cite[Theorem 7.1]{Y}
    The coarse Baum--Connes conjecture holds for proper metric spaces with finite asymptotic dimension.
\end{Theorem}

Therefore, combining \cref{thm:4Thurstongeom-finite-asdim}, \cref{thm:3Alex-asdim-bound-3}, and \cref{thm:CBC} we obtain a proof of the following result.

\begin{Corollary}\label{cor:Coarse-BC-conj}
\begin{enumerate} [i)]
    \item Let $X$ be an oriented closed $4$-manifold that is either geometric or admits a geometric decomposition as in \cref{thm:4Thurstongeom-finite-asdim}. Then the coarse Baum--Connes conjecture holds for $\pi_1(X)$.
    \item Let $Y$ be a closed $3$-dimensional Alexandrov space. Then the coarse Baum--Connes conjecture holds for $\pi_1(Y)$.
\end{enumerate}    
\end{Corollary}

\subsection{Novikov conjecture}\label{subsec:Novikov}

Let $M$ be a manifold and $\Gamma=\pi_1(M)$. If $M$ is oriented, a rational cohomology class $x\in H^{\ast}(B\Gamma, \mathbf{Q})$ defines a rational characteristic number, called a higher signature (see \cite{Y}):
\begin{equation*}
    \sigma_{x}(M,u) = \ \langle \mathcal{L}(M)\cup u^{\ast} (x), [M]\rangle \in \mathbf{Q}
\end{equation*}
Here $\mathcal{L}(M)$ is the Hirzebruch $\mathcal{L}$-genus, and $u: M \longrightarrow B\Gamma$ the classifying map.
The Novikov conjecture posits that all higher signatures are invariants of oriented homotopy equivalences over $B\Gamma$.

As a way of explaining how a good picture of a metric space could be "drawn" inside a Hilbert space, Gromov \cite{Gro2} introduced the following concept.

\begin{Definition}
    Let $(H,d_{H})$ be a Hilbert space and $(X, d)$ a metric space. 
    A map $f:X\to H$ is a coarse embedding into $H$ if there exist non-decreasing functions $\rho_1$ and $\rho_2$ on $[0, \infty )$ such that:
    \begin{enumerate}
        \item $\rho_1(d(x,y)) \leq d_{H}(f(x), f(y)) \leq \rho_1(d(x,y))$, for all $x,y$ in X.
        \item $\lim_{r\to +\infty} \rho_1(r) =  +\infty$.    
    \end{enumerate}
\end{Definition}

Coarse embeddability of a countable group into a Hilbert space is independent of the choice of proper length metrics \cite{Gro2}. 
Crucially for us, groups with finite asymptotic dimension are coarsely embeddable into a Hilbert space \cite{Yu2}.
Moreover, the next result follows from Yu \cite{Yu2}, Higson \cite{Higson}, and Skandalis--Tu--Yu \cite{STY} (see \cite[\S 3]{Yu}).

\begin{Theorem}\label{thm:HilbertEmbedd-implies-Novikov}
    The Novikov conjecture holds if the fundamental group of a manifold is coarsely embeddable into a Hilbert space.
\end{Theorem}

Therefore, combining \cref{thm:4Thurstongeom-finite-asdim}, \cref{thm:3Alex-asdim-bound-3}, and \cref{thm:HilbertEmbedd-implies-Novikov}, we obtain the following result. 

\begin{Corollary}\label{cor:Novikov-Conj}
\begin{enumerate} [i)]
    \item Let $X$ be an oriented closed $4$-manifold that is either geometric or admits a geometric decomposition as in \cref{thm:4Thurstongeom-finite-asdim}. Then the Novikov conjecture holds for $X$.
    \item Let $Y$ be a closed $3$-dimensional Alexandrov space. Then the Novikov conjecture holds for $Y$.
    \end{enumerate}  
\end{Corollary}

To the best of our knowledge, the second item in \cref{cor:Novikov-Conj} above is new.

\subsection{Zero in the spectrum}\label{subsec:Zero}

Recall that the Laplace--Beltrami operator $\Delta_{p}$, with $0\leq p \leq n$, of a complete oriented Riemannian $n$-manifold acts on square-integrable forms. 
It is an essentially self-adjoint positive operator, so its spectrum is a subset of the positive reals. 
A space $X$ is said to be uniformly contractible if, for each $R>0,$ there exists some $S>R$ such that, for all $x\in X$, the ball $B(x,R)$ is contractible within $B(x,S)$. 
Gromov's zero in spectrum conjecture asks if the spectrum of $\Delta_{p}$  of a uniformly contractible Riemannian $n$-manifold contains zero, for any $0\leq p \leq n$ (see \cite{Lott}). 
As a consequence of \cref{thm:CBC}, Yu showed: 

\begin{Corollary}\label{cor:Yu-zero-spec} \cite[Corollary 7.4]{Y}
        Gromov's zero-in-the-spectrum conjecture holds for uniformly contractible Riemannian manifolds with finite asymptotic dimension.
\end{Corollary}

Aspherical manifolds have uniformly contractible universal covering spaces, therefore the following holds:

 \begin{Corollary}\label{cor:zerospec}
 Let $Z$ be either an aspherical manifold from \cref{thm:4Thurstongeom-finite-asdim}, or a closed aspherical $3$--dimensional Alexandrov space. Then, there exists a $p\geq 0$, such that zero belongs to the spectrum of the Laplace--Beltrami operator $\Delta_{p}$ acting on square-integrable $p$-forms of the universal cover $\widetilde{Z}$ of $Z$.
 \end{Corollary}

 \begin{proof}
     Observe that, by \cref{cor:Yu-zero-spec},  the result holds for an aspherical manifold $X$ from  \cref{thm:4Thurstongeom-finite-asdim}.
     
     Let $Y$ be a closed, aspherical, $3$--dimensional Alexandrov space. 
     Then, by the universal property in \cref{lem:3Alex-Univ-Cover}, its universal covering space $\widetilde{Y}$ is also the universal cover of the manifold $Y^{*}$ that is a (potential) double branched cover of $Y$. So $\widetilde{Y}$ is a manifold.
     Moreover, we already know that the asymptotic dimension of $\widetilde{Y}$ is finite, because it equals that of $\pi_{1}(Y) \leq 3$.
     Therefore, again by \cref{cor:Yu-zero-spec}, Gromov's zero-in-the-spectrum conjecture holds for $\widetilde{Y}$. 
 \end{proof}

\subsection{Yamabe invariant}\label{subsec:Yamabe}

Obtaining bounds, or exact computations, of the Yamabe invariant is a notoriously difficult problem.
Schoen \cite{Sc} showed that:
\begin{enumerate}[i)]
    \item $M$ has $\mathcal{Y}(M)>0$ if and only if it admits a positive scalar curvature smooth metric.
    \item If $M$ admits a volume collapsing sequence of metrics with bounded curvature, then $\mathcal{Y}(X)\geq 0$. 
\end{enumerate}
A notable result by Petean states that every simply connected smooth compact manifold of dimension greater than four has non-negative Yamabe invariant \cite{Petean}.

As previously mentioned, Yu showed that an aspherical manifold with fundamental group of finite asymptotic dimension does not admit a metric of positive scalar curvature \cite{Y}. 

We will now recall a notion first introduced by Gromov \cite{Gro}, that generalizes the effect of having a circle action, in terms of vanishing of various invariants of smooth manifolds (cf. \cite{PaternainPetean}).

An $\mathcal{F}$-structure on a closed manifold $M$ is given by,
\begin{enumerate}
\item{ A finite open cover $\{ U_1, ..., U_{N} \} $;}
\item{ $\pi_{i}\co\widetilde{U_{i}}\rightarrow U_{i}$ a finite Galois covering with group of deck transformations $\Gamma_{i}$, $1\leq i \leq N$;}

\item{ A smooth torus action with finite kernel of the $k_{i}$-dimensional torus, \\ $\phi_{i}\co T^{k_{i}}\rightarrow {\rm{Diff}}(\widetilde{U_{i}})$, $1\leq i \leq N$;}

\item{ A homomorphism $\Psi_{i}\co \Gamma_{i}\rightarrow {\rm{Aut}}(T^{k_{i}})$ such that
\[ \gamma(\phi_{i}(t)(x))=\phi_{i}(\Psi_{i}(\gamma)(t))(\gamma x) \]
for all $\gamma \in \Gamma_{i}$, $t \in T^{k_{i}}$ and $x \in \widetilde{U_{i}}$; }

\item{ For any finite sub-collection $\{ U_{i_{1}}, ..., U_{i_{l}} \} $ such that  $U_{i_{1}\ldots i_{l}}:=U_{i_{1}}\cap \ldots \cap U_{i_{l}}\neq\emptyset$ the following compatibility condition holds: let $\widetilde{U}_{i_{1}\ldots i_{l}}$ be the set of points $(x_{i_{1}}, \ldots , x_{i_{l}})\in \widetilde{U}_{i_{1}}\times \ldots \times \widetilde{U}_{i_{l}}$ such that $\pi_{i_{1}}(x_{i_{1}})=\ldots = \pi_{i_{l}}(x_{i_{l}})$. The set $\widetilde{U}_{i_{1}\ldots i_{l}}$ covers $\pi_{i_{j}}^{-1}(U_{i_{1}\ldots i_{l}}) \subset \widetilde{U}_{i_{j}}$ for all $1\leq j \leq l$. Then we require that $\phi_{i_{j}}$ leaves $\pi_{i_{j}}^{-1}(U_{i_{1}\ldots i_{l}})$ invariant and it lifts to an action on $\widetilde{U}_{i_{1}\ldots i_{l}}$ such that all lifted actions commute. }
\end{enumerate}

The second named author showed that manifolds in \cref{cor:zero-Yamabe} admit an $\mathcal{F}$-structure.

\begin{Theorem}[Theorems A \& B \cite{SS09}]\label{thm:4Dgeomz-Fstructure}
    Let $X$ be a manifold that is either geometric or that it admits a geometric decomposition into pieces modelled on one of the following geometries:
\[
    \begin{array}{cccccccc}
     \mathbb{S}^{4}, & \mathbb{C}{\rm P}^{2}, & \mathbb{S}^{3}\times \mathbb{E}, & \mathbb{H}^{3}\times \mathbb{E},  & 
 \Sl \times \mathbb{E}, & \mathbb{N}il^{3}\times \mathbb{E},  & \mathbb{N}il^{4}, &  \mathbb{S}ol^{4}_{1}, \\
 \mathbb{S}^{2} \times\mathbb{E}^{2}, & \mathbb{H}^{2}\times\mathbb{E}^{2},
& \mathbb{S}ol^{4}_{m,n}, & \mathbb{S}ol^{4}_{0},  &
 \mathbb{S}^{2}\times \mathbb{S}^{2},  & \mathbb{S}^{2}\times\mathbb{H}^{2}, & \mathbb{E}^{4}, & \mathbb{F}^4
\end{array} 
\]
Then $X$ admits an $\mathcal{F}$-structure.
\end{Theorem}

The connection between the existence of $\mathcal{F}$-structures and bounds for the Yamabe invariant is given by the next theorem of Paternain and Petean.

\begin{Theorem}[Theorem 7.2 \cite{PaternainPetean}]\label{thm:Fstruct-YamabeZero}
If a closed smooth manifold $X$ admits an $\mathcal{F}$-structure, $\dim\, X> 2$, then $\mathcal{Y}(X)\geq 0$.
\end{Theorem}

We are now ready to present a proof of \cref{cor:zero-Yamabe}.

\begin{proof}[Proof of \cref{cor:zero-Yamabe}]
    By  \cref{thm:4Dgeomz-Fstructure} the manifolds with geometric pieces modeled on these geometries admit an $\mathcal{F}$-structure.    
    Then, by Paternain and Petean's \cref{thm:Fstruct-YamabeZero}, its Yamabe invariant is non-negative.
    From \cref{thm:4Thurstongeom-finite-asdim} and Yu's celebrated result \cite{Y} it follows that such a manifold $X$ does not admit a metric of positive scalar curvature.  
    Then Schoen's result mentioned above implies $\mathcal{Y}(X)\leq 0$.
    Therefore $\mathcal{Y}(X)=0$.    

    Now consider the connected sums of $X$ with $S^3\times S^1$. 
    We appeal to a result of Petean, who showed that performing $0$-dimensional surgery on $X$ leaves the Yamabe invariant unchanged \cite[Proposition 3]{Petean1}. 
    Iterating this last argument yields the result for any finite number of  connected sums with $S^3\times S^1$, as claimed. 
\end{proof}

We now include a proof of \cref{lem:ComplexSurfaces-ZeroYam}.

\begin{proof}[Proof of \cref{lem:ComplexSurfaces-ZeroYam}]
By the work of Paternain and Petean on collapsing of compact complex surfaces, $X$ admits an $\mathcal{F}$-structure \cite[Theorems A \& B]{PaternainPetean2}. 
Hence, \cref{thm:Fstruct-YamabeZero} yields $\mathcal{Y}(X)\geq 0$.
Now, by Chodosh and Li \cite{CL23}, and Gromov \cite{Gro3}, $X$ does not admit a metric of positive scalar curvature.
Thus, by the previously mentioned results we obtain $\mathcal{Y}(X)\leq 0$.
Therefore we conclude $\mathcal{Y}(X)= 0$.
\end{proof}

\section{Proofs of our main results}\label{sec:main-proofs}

We start with the following lemma, needed for the proofs of our main results. 

 \begin{Lemma}\label{lem:asdim-3D-geometric}
       Let $Y$ be a compact $3$-manifold that is geometric in the sense of Thurston, then ${\rm asdim}\ \pi(Y) \leq 3$.
    \end{Lemma}

    \begin{proof}
    This can be verified for each of the model geometries, which we group as follows:
    \begin{enumerate}
         \item $\mathbb{E}^{3}$\ ; this case follows from \cref{lem:finite-asdim}, item (i).
         \item $\mathbb{H}^{3}$\ ; this case follows from \cref{lem:finite-asdim}, item (ii).
         \item $\mathbb{S}^{3}$\ ; these groups are finite, it follows from \cref{lem:finite-asdim}, item (iii).
        \item $Nil^3, Sol^3, \Sl$\ ; these cases are covered by \cref{lem:asdim-LieLattice}. 
        The geometries $\mathbb{E}^{3}, Nil^3$ and $Sol^3$ are all Lie groups. Notice that $\Sl$ is the universal cover of the $3$-dimensional Lie group ${\rm SL}_{2}$ of all $2\times 2$ matrices with determinant $1$. 
        As $\Sl$ admits an invariant metric under left or right multiplication, then $\Sl$ is also a Lie group.

        \item $\mathbb{S}^{2}\times \mathbb{E}, \mathbb{H}^{2}\times \mathbb{E}$\ ; for these geometries the proof follows from the previously mentioned result for $\mathbb{S}^{2}$ and $\mathbb{H}^{2}$ in \cref{lem:asdim-2D} in combination with the bound for products of spaces found in \cref{prop:asdim-products}.        
    \end{enumerate}
    Therefore, in all the possible cases we obtain that ${\rm asdim}\ \pi(Y) \leq 3$.
    \end{proof}

\subsection{Proof of \cref{thm:4Thurstongeom-finite-asdim}}

\begin{proof}[Proof of Theorem \ref{thm:4Thurstongeom-finite-asdim}]

First we will prove the upper bound for ${\rm asdim}\pi_1(M)$ in both cases, the geometric manifolds and the manifolds with geometric decomposition. 

{\it 1) Geometric Manifolds.} 

\indent {\it 1.1) Finite fundamental groups; $  \mathbb{S}^{4},  \mathbb{S}^{2}\times \mathbb{S}^{2}, \mathbb{C}{\rm P}^{2}$.}
Let $\Gamma$ be a finite group, then $\Gamma$ is finitely generated. Hence, by item (3) in Lemma \ref{lem:finite-asdim} the asymptotic dimension ${\rm asdim}\, \Gamma =0$. 
Therefore, as the fundamental groups of geometric manifolds modelled on $ \mathbb{S}^{4},  \mathbb{S}^{2}\times \mathbb{S}^{2},$ or $ \mathbb{C}{\rm P}^{2}$ are finite, they have asymptotic dimension zero.\\

{\it 1.2) Quotients of Lie groups.}
 The asymptotic dimension of quotients of simply connected Lie groups can be effectively bounded.
By Lemma \ref{lem:asdim-LieLattice} a cocompact lattice $\Gamma$ in a connected Lie group $G$ with maximal compact subgroup $K$ satisfies ${\rm asdim}\, \Gamma = \dim (G/K)$.
Hence we obtain that ${\rm asdim}\ \pi_1(X)\leq 4$ for geometric manifolds $X$ modelled on the geometries 
$ \mathbb{N}il^{3}\times \mathbb{E},\   \mathbb{N}il^{4},\   \mathbb{S}ol^{4}_{1},\  \mathbb{S}ol^{4}_{m,n},\  \mathbb{S}ol^{4}_{0}$, or $ \mathbb{E}^{4} $.\\

{\it 1.3) Product geometries; $ \mathbb{S}^{3}\times \mathbb{E},\  \mathbb{H}^{3}\times \mathbb{E},\ \Sl \times \mathbb{E},\ \mathbb{S}^{2} \times\mathbb{E}^{2},\  \mathbb{H}^{2}\times\mathbb{E}^{2},\  \mathbb{S}^{2}\times\mathbb{H}^{2}$, $\mathbb{H}^{2}\times\mathbb{H}^{2}$.}
By Proposition \ref{prop:asdim-products}, we know that the asymptotic dimension of a product of coarse spaces 
is bounded by the sum of the asymptotic dimension of each space. 
Therefore, for all the product geometries $ \mathbb{S}^{3}\times \mathbb{E},\  \mathbb{H}^{3}\times \mathbb{E},\ \mathbb{S}^{2} \times\mathbb{E}^{2},\  \mathbb{H}^{2}\times\mathbb{E}^{2},\  \mathbb{S}^{2}\times\mathbb{H}^{2}$ and $\mathbb{H}^{2}\times\mathbb{H}^{2}$, we have that their asymptotic dimension is bounded above by the sum of the asymptotic dimension of their factors. 
Therefore, by \cref{lem:asdim-2D}, \cref{lem:asdim-3D-geometric}, and \cref{prop:asdim-products} the asymptotic dimension of each of these product geometries is at most $4$.\\

{\it 1.4) $\mathbb{H}^4$ and $\mathbb{H}^2_{{\bf C}}$ manifolds;}
Observe that compact $\mathbb{H}^4$ or $\mathbb{H}^2_{{\bf C}}$ manifolds have hyperbolic fundamental groups, so by \cref{thm:hyp-fin-asdim} above their asymptotic dimension is finite. 
Finite volume manifolds modeled on these geometries, truncated to be used as pieces of a geometric decomposition, are relatively hyperbolic with respect to their peripheral structure, that is, the systems of fundamental groups of their boundary components. 
Such groups have finite asymptotic dimension by \cref{thm:rel-hyp-fin-asdim}. 
Moreover, we showed in \cref{cor:Nagata-complexhyperbolic} that complex hyperbolic pieces have Nagata dimension 4. That real hyperbolic pieces have asymptotic dimension 4 follows from item {\it ii)} in  \cref{lem:finite-asdim}.\\

{\it 1.5) $\mathbb{F}^4$ manifolds;}
For the case of $\mathbb{F}^4$, the extension result of Bell--Dranishnikov in \cref{thm:extension}, applied to a short exact sequence of the fundamental group, yields the desired bound.
Let $X$ be a manifold modeled on $\mathbb{F}^4$, then $\pi_1(X)$ is isomorphic to a lattice in $\mathbf{R}^2\ltimes {\rm SL}(2,\mathbf{R})$. 
Let $\overline{\pi_1(X)}$ be the image of $\pi_1(X) $ in ${\rm SL}(2,\mathbf{R})$.
Observe that $\overline{\pi_1(X)}= \pi_1(X) / \mathbf{Z}^2$. 
This exhibits $X = \mathbb{F}^4 / \pi_1(X)$ as an elliptic surface, over the base $B=\mathbb{H}^2 / \overline{\pi_1(X)}$, where $B$ is a non-compact orbifold \cite[p.150]{Wall}.

 The identity component of ${\rm Iso}(\mathbb F^4)$  is the semidirect product
 $\mathbf{R}^{2}\ltimes_{\alpha}{\rm SL}(2,\mathbf{R})$, with $\alpha$ the natural action of
 ${\rm SL}(2,\mathbf{R})$ on $\mathbf{R}^{2}$. 
Let $p \co \mathbf{R}^{2}\ltimes_{\alpha}{\rm SL}(2,\mathbf{R})\to {\rm SL}(2,\mathbf R)$
be the projection homomorphism. 
 The manifold $X$ is diffeomorphic to tue quotient of $T^2\times\mathbb H^2$ under the action of $p(\pi_1(X))$, acting on $ T^2$ through $\psi \co p(\pi_1(X))\to {\rm Aut}( T^2=\mathbb R^2/(\pi_1(X)\cap\mathbb R^2))$, and on $\mathbb H^2$ in the usual way. 

The quotient $B:=\mathbb H^2/p(\pi_1(X))$
is a finite volume hyperbolic $2$-orbifold, hence $X$ is an orbifold bundle over $B$.
If $B$ is smooth, i.e. $p(\pi_1(X))$ acts without fixed points, then $M$ is a torus bundle over $B$ with structure group ${\rm SL}(2,{\bf Z})$ and $\psi$ is precisely its holonomy.

In general, the manifold $X$ is a $T^2$ fibration over a noncompact, finite area, hyperbolic orbifold $B$ ( see \cite[p.150]{Wall} and \cite[\S 10.1]{ SS09}). 
Then, the fundamental group $\pi_1(X)$ can be written as an extension, $ \pi_1(T^2)\to \pi_1(X)\to \pi_1(B)$.
Thus \cref{eq:extension} of the extension Theorem \ref{thm:extension} implies 
${\rm asdim}\ \pi_1(X) \leq {\rm asdim}\ \pi_1(T^2) + {\rm asdim}\ \pi_1(B) \leq 2+2.$

Therefore, our arguments have now covered all the possible cases and we conclude that all the $4$-dimensional geometric manifolds have asymptotic dimension at most $4$.\\

 {\it 2) Manifolds with a geometric decomposition.}

As we explained in \cref{subsec:graph-gps} above, the fundamental group $\pi_1(X)$ of a geometrizable $4$-manifold $X$ with a proper and $\pi_1$-injective geometric decomposition is  isomorphic to a graph of groups. 
By \cref{thm:amalgams-fin-asdim} the asymptotic dimension of a graph of groups is finite provided each vertex group has finite asymptotic dimension, moreover we have computed the explicit bound we need. 

We will now cover each of the possible geometric decompositions, in the same order of Hillman's \cref{thm:aspherical-manifolds} above.\\

{\it 2.1) $X$ is the total space of an orbifold bundle with general fiber $S^2$ over a hyperbolic 2-orbifold;}
Notice that, by \cref{thm:proper-action}, the relevant fiber and base orbifold groups have asymptotic dimension at most $2$.
Consider the decomposition of $X$ into its geometric pieces $X_{i}$, $i \in \{ 1, \ldots , k\}$ (see \cite{Hil2}).
Then, the arguments explained above for the geometric cases yield ${\rm asdim}\ \pi_1(X_{i})\leq 4$.
As $X$ is compact $k< \infty$, so the finite union finite union \cref{thm:finite-union} implies $ {\rm asdim}\ \pi_1(X) \leq 4$.\\

{\it 2.2) Manifolds that decompose into  $\mathbb{H}^{2}\times\mathbb{H}^{2}$ pieces;}

These manifolds have two kinds of decompositions, called irreducible if the boundary inclusions into each piece is $\pi_1$-injective, or reducible otherwise.

In the irreducible case we obtain a decomposition of the fundamental group into a graph of groups, and the result follows as in other similar cases.

In the reducible case we use the finite union \cref{thm:finite-union}. 
Apply it first to a couple of contiguously glued $\mathbb{H}^{2}\times\mathbb{H}^{2}$-pieces.
Then perform induction over the number of pieces of the geometric decomposition, to obtain the desired upper bound.
Hence in both the reducible and irreducible cases we have that $ {\rm asdim}\ \pi_1(X) \leq 4$.\\

{\it 2.3) Manifolds that decompose into $\mathbb{H}^{4},\ \mathbb{H}^{3}\times \mathbb{E}^{1}, \ \mathbb{H}^{2}\times \mathbb{E}^{2}$, or $\Sl \times \mathbb{E}^{2}$ pieces;}  
First, observe that a manifolds $X$ that decompose into pieces modelled on the hyperbolic geometry $\mathbb{H}^{4}$ are relatively hyperbolic. 
Their ends are either flat or nilpotent, and therefore the fundamental groups of each geometric piece has finite asymptotic dimension by Theorem \ref{thm:rel-hyp-fin-asdim}, because they are relatively hyperbolic. 
Notice that the fundamental groups $\pi_1(Y)$ of flat or nilpotent $3$-manifolds $Y$ have ${\rm asdim}\ \pi_1(Y) = 3$, by Lemma \ref{lem:asdim-LieLattice}.

Now, for the case that the manifold $X$ decompose into pieces modelled on the geometries $\ \mathbb{H}^{3}\times \mathbb{E}^{1}, \ \mathbb{H}^{2}\times \mathbb{E}^{2}$ or $\Sl \times \mathbb{E}^{2}$ we already proved that the asymptotic dimension of all this product geometries is at most $4$. 
Therefore, in both cases by Theorem \ref{thm:amalgams-fin-asdim} the asymptotic dimension of $\pi_1(X)$ is bounded above by $4$. \\

{\it 2.4) Manifolds that decompose into $\mathbb{H}_{\C}^{2}$ or $\mathbb{F}^{4}$ pieces;}  
If a manifold $X$ decompose into pieces modelled on the hyperbolic geometry $\mathbb{H}_{\C}^{2}$, then again it is relatively hyperbolic. 
So, as in the case of $\mathbb{H}^{4}$, the fundamental group of each geometric piece has finite asymptotic dimension. 

For the case that $X$ decompose into pieces modelled on the geometry $\mathbb{F}^{4}$ we know that the asymptotic dimension of each piece is bounded by above by $4$.
By Theorem \ref{thm:amalgams-fin-asdim}, the asymptotic dimension of $\pi_1(X)$ is bounded above by 4.\\

Therefore, we have shown ${\rm asdim}\ \pi_1(M)\leq 4$ when $X$ is a closed orientable $4$-manifold that is geometric or admits a geometric decomposition in the sense of Thurston. \\

{\it 3) Equality for aspherical manifolds.}
Now we will prove the lower bound ${\rm asdim}\ \pi_1(M)\geq 4$ in the case of aspherical manifolds, which will imply the equality we claim.
By Theorem \ref{thm:aspherical-manifolds}, we know that a geometric manifold modelled on $\mathbb{H}^{3}\times \mathbb{E},\ \mathbb{H}^{2}\times \mathbb{E}^{2}, \ \mathbb{H}^{4}, \ \mathbb{H}^{2}\times \mathbb{H}^{2}, \ \mathbb{H}^{2}_{\C}, \Sl \times \mathbb{E}, \ \mathbb{N}il^{3}\times \mathbb{E},  \ \mathbb{N}il^{4}, \  \mathbb{S}ol^{4}_{1}, \ \mathbb{S}ol^{4}_{m,n},$ or $ \ \mathbb{S}ol^{4}_{0}$ is aspherical. 
Hence, Lemma \ref{lmm:asph-cohdim} implies that the asymptotic dimension of such a fundamental group is bounded below by its cohomological dimension. 
Observe that the cohomological dimension of $\pi_1(M)$ is equal to the dimension of $M$, so 
${\rm asdim}\ \pi_1(M) \geq \dim M = 4$.
This concludes the proof for geometric manifolds.

In the cases of manifolds with a geometric decomposition, by items (4) and (5) of Theorem \ref{thm:aspherical-manifolds}, we know that if the pieces of the geometric decomposition have geometries $\mathbb{H}^{4},\ \mathbb{H}^{3}\times \mathbb{E}^{1}, \ \mathbb{H}^{2}\times \mathbb{E}^{2}$, $\widetilde{\mathbb{SL}} \times \mathbb{E}^{2}$, $\mathbb{H}_{\C}^{2}$ or $\mathbb{F}^{4}$, then the manifold is aspherical.
Using Lemma \ref{lmm:asph-cohdim} again, we obtain that the lower bound for their asymptotic dimension is $4$.

Therefore the equality follows for both cases, covering all the possible aspherical manifolds as we claimed.

Finally, we mention the effect of connected sums on the asymptotic dimension.
Observe that taking connected sums of manifolds corresponds to performing free products at the level of fundamental groups. 
So that if the pieces of the connected sum have finite asymptotic dimension, then the resulting connected sum also has finite asymptotic dimension. 
Moreover, the upper bound in this case remains the same, according to Theorem \ref{thm:amalgams-fin-asdim}.
This concludes our proof. \end{proof}

\subsection{Proof of \cref{thm:3mfds-asdim-leq-3}}

We will now present a proof of \cref{thm:3mfds-asdim-leq-3}.

\begin{proof}[Proof of \cref{thm:3mfds-asdim-leq-3}]
    The success of Thurston's geometrization program implies that $\pi_1(Y)$ may be presented as a graph of groups $\mathcal{G}_{Y}$; each vertex group $V_{i}$ is a discrete group of isometries of one of the eight model geometries, while the edge groups $E_{i,j}$---between the vertices $V_{i}$ and $V_{j}$---are surface groups.
    By \cref{lem:asdim-3D-geometric} the asymptotic dimension of the groups $V_{i}$ is bounded above by $3$, that is, ${\rm asdim}\ V_{i} \leq 3$.

     Continuing with the proof, observe that a finite graph of groups is isomorphic to an interated amalgamated product. 
 Therefore, by \cref{thm:amalgams-fin-asdim}, and lemmata \ref{lem:asdim-3D-geometric} and \ref{lem:asdim-2D}, we obtain
 ${\rm asdim}\ \mathcal{G}_{Y} \leq \max \{ {\rm asdim}\ V_{i} ,\  {\rm asdim}\ E_{i,j} +1\} \leq 3$.
 In the non-orientable case, consider the orientation double cover to obtain the same result.

Finally, in the aspherical case, Lemma \ref{lmm:asph-cohdim} yields $3$ as the lower bound for the asymptotic dimension.
 
\end{proof}

\subsection{Proof of \cref{thm:3Alex-asdim-bound-3}}

Now  we continue with a proof of \cref{thm:3Alex-asdim-bound-3}.

\begin{proof}[Proof of \cref{thm:3Alex-asdim-bound-3}]
  
Let $Y$ be a compact $3$-dimensional Alexandrov space.
Then, as explained above, by \cite{GGG}, there exist both a smooth Riemannian $3$-manifold $Y^{*}$, and an isometric involution $\iota$ of $Y^{*}$, such that $Y\simeq Y^{*} / \{ y \cong \iota(y) \} $, with $y\in Y^{*}$.
Write $\widetilde{Y}$ for the universal covering of $Y^{*}$.
Then the fundamental group $\pi_1 (Y^{*})$, seen as a group of deck transformations, acts on $\widetilde{Y}$.
Moreover, observe that, as $\iota$ is an isometric involution acting on $Y^{*}$, it lifts to an action on $\widetilde{Y}$.

Denote by $\Gamma$ the group formed by composing the action of $\pi_1 Y^{*}$ on $\widetilde{Y}$, with the action of $\iota$ on $Y^{*}$.
The orbit equivalence classes of $\Gamma$, acting on $\widetilde{Y}$, present $Y$ as a quotient space.
Recall that, by \cref{lem:3Alex-Univ-Cover}, $\widetilde{Y}$ is the unique universal cover of both $Y^{*}$ and $Y$.
Therefore $\Gamma$ is isomorphic to the fundamental group of $Y$.

We claim that $\Gamma$ acts properly on $\widetilde{Y}$, because the isotropy groups are of two kinds only:
\begin{enumerate}[i)]
    \item The trivial group, for the non-singular points of $Y$, whose space of directions is homeomorphic to a ball.
    \item Isomorphic to ${\bf Z} / 2$, for the singular points $\mathcal{S}(Y)$, whose space of directions is homeomorphic to a projective plane.
\end{enumerate}
As these two cases cover all possible types of isotropy groups, the action is proper, as claimed. 
Hence, the group acts properly, and also isometrically, on the proper metric space $\widetilde{Y}$. Therefore \cref{thm:proper-action} implies ${\rm asdim}\ \Gamma \leq {\rm asdim}\ \widetilde{Y}$, and by \cref{eqn:universal-asdim} and \cref{thm:3mfds-asdim-leq-3}, ${\rm asdim}\ \widetilde{Y} \leq 3$. 
\end{proof}

\providecommand{\bysame}{\leavevmode\hbox to3em{\hrulefill}\thinspace}

\end{document}